\documentclass{amsart}

\usepackage{cleveref}
\usepackage{enumitem}
\usepackage{amsfonts, amssymb}
\usepackage{tikz}
\usetikzlibrary{decorations.markings}
\usepackage{mathtools}
 
 \tikzset{->-/.style={decoration={
  markings,
  mark=at position .5 with {\arrow{>}}},postaction={decorate}}}

\newtheorem{thm}{Theorem}[section] 
\newtheorem{lem}[thm]{Lemma}
\newtheorem{cor}[thm]{Corollary}
\newtheorem{prop}[thm]{Proposition}
\newtheorem{alg}[thm]{Algorithm}

\newcommand{\zz}{\mathbb{Z}}
\newcommand{\rr}{\mathbb{R}}
\newcommand{\nn}{\mathbb{N}}
\newcommand{\cc}{\mathbb{C}}
\newcommand{\qq}{\mathbb{Q}}
\newcommand{\hh}{\mathbb{H}}
\newcommand{\oo}{\mathcal{O}}
\newcommand{\pp}{\mathcal{P}}

\title[Discrete and free two-generated subgroups of ${\rm SL_2}$]{Discrete and free two-generated subgroups of ${\rm SL_2}$ over non-archimedean local fields}
\author{Matthew J. Conder}
\address{M. J. Conder, Department of Pure Mathematics and Mathematical Statistics, Centre for Mathematical Sciences, University of Cambridge, Wilberforce Road, Cambridge, CB3 0WB, United Kingdom}
\email{mjc271@cam.ac.uk}
\address{F. Paulin, Laboratoire de math\'ematique d'Orsay, UMR 8628 Univ. Paris-Sud et CNRS, Universit\'e Paris-Saclay, 91405 Orsay Cedex, France}
\email{frederic.paulin@math.u-psud.fr}

\begin{document}

\begin{abstract}
We present a practical algorithm which, given a non-archimedean local field $K$ and any two elements $A,B\in {\rm SL_2}(K)$, determines after finitely many steps whether or not the subgroup $\langle A, B \rangle\le {\rm SL_2}(K)$ is discrete and free of rank two. This makes use of the Ping Pong Lemma applied to the action of ${\rm SL_2}(K)$ by isometries on its Bruhat-Tits tree. The algorithm itself can also be used for two-generated subgroups of the isometry group of any locally finite simplicial tree, and has applications to the constructive membership problem. In an appendix joint with Fr\'ed\'eric Paulin, we give an erratum to \cite[Proposition 1.6]{Paulin}, which details some translation length formulae that are fundamental to the algorithm.
\end{abstract}

\maketitle

\section{Introduction}

The problem of deciding whether or not two elements of ${\rm SL_2}(\rr)$ generate a free group of rank two has been widely studied in the literature. For instance, the subgroups generated by matrices of the form $\left[ \begin{array}{cc}
1 & \alpha \\
0 & 1 \end{array} \right]$,
$\left[ \begin{array}{cc}
1 & 0 \\
\alpha & 1 \end{array} \right]$ are known to be free of rank two whenever $|\alpha |\ge 2$; this is an easy consequence of the Ping Pong Lemma, applied to the action of ${\rm SL_2}(\rr)$ on the hyperbolic plane $\hh^2$ via M\"{o}bius transformations. On the other hand, there are many rational values of $\alpha$ in the interval $(-2,2)$ for which the subgroup generated by the above matrices is not free, and it is an open question to decide whether or not this holds for every such rational $\alpha$; see, amongst other papers, \cite{B} and \cite{LU}.

A key observation in \cite{N} is that arguments involving the Ping Pong Lemma can show that some two-generated subgroups of ${\rm SL_2}(\rr)$ which are free are also discrete, with respect to the topology inherited from $\rr^4$. This helped lead to the discovery of necessary and sufficient conditions, depending on matrix trace, for a two-generated subgroup of ${\rm SL_2}(\rr)$ (or, equivalently, of ${\rm PSL_2}(\rr)$) to be discrete and free of rank two; see \cite{P} or \cite {R2}. Moreover, given any two elements $A,B \in {\rm SL_2}(\rr)$, Nielsen transformations can be performed in a `trace minimising' manner to determine whether or not these conditions are satisfied for the subgroup $\langle A, B \rangle \le {\rm SL_2}(\rr)$. This observation (also made in \cite{G} in the context of determining discreteness) forms the basis of a practical algorithm given explicitly in \cite{SL2}, which determines after finitely many steps whether or not a given two-generated subgroup of ${\rm SL_2}(\rr)$ (or ${\rm PSL_2}(\rr)$) is discrete and free. It is also noted in \cite{SL2} that this algorithm can be used to solve the constructive membership problem for discrete and free two-generated subgroups of ${\rm SL_2}(\rr)$ or ${\rm PSL_2}(\rr)$; namely, given such a subgroup $G$ and an element $X$ in the corresponding overgroup, one can determine algorithmically whether or not $X$ lies in $G$, and if it does, give an explicit representation of $X$ as a word in the generators of $G$.

Discrete and free two-generated subgroups of ${\rm SL_2}$ over other fields, particularly other locally compact fields, are not as well studied. There has been some work done in the case of ${\rm SL_2}(\cc)$ (for instance, see \cite{Bow}) but the action of this group on hyperbolic space $\hh^3$ is much more complicated to study. Over a non-archimedean local field $K$, however, the group ${\rm SL_2}(K)$ acts by isometries and without inversions on the corresponding Bruhat-Tits tree, and such actions on simplicial trees are very well understood. Given two elements $A,B \in {\rm SL_2}(K)$, we will show that Nielsen transformations can be performed in a `translation length minimising' manner until either the subgroup $\langle A,B \rangle \le {\rm SL_2}(K)$ is shown to contain an elliptic element (which is either of finite order or generates an indiscrete infinite cyclic subgroup), or hyperbolic generators of $\langle A,B \rangle$ are found which satisfy the hypotheses of the Ping Pong Lemma. This helped us form the basis of a practical algorithm (\Cref{treealg}) which determines after finitely many steps whether or not a given two-generated subgroup of ${\rm SL_2}(K)$ (or, equivalently, of ${\rm PSL_2}(K)$) is discrete and free. We will show that this algorithm can also be used more generally (in the context of isometry groups of locally finite simplicial trees) and gives a further algorithm solving the constructive membership problem for such groups that are discrete and free.

\medskip

In Section 2, we provide some background information on non-archimedean local fields and the group ${\rm SL_2}(K)$ defined over such a field $K$. We describe the Bruhat-Tits tree associated to such groups and some general theory of groups acting on simplicial trees by isometries and without inversions.

Section 3 details the key results leading to \Cref{treealg}; in particular, we show that a discrete and free subgroup of ${\rm SL_2}(K)$ cannot contain any elliptic elements, and present a form of the Ping Pong Lemma that gives conditions for a pair of hyperbolic elements to generate a discrete and free subgroup. We also give some important translation length formulae, one of which corrects a formula given in \cite[Proposition 1.6]{Paulin}. In the appendix, joint with the author of \cite{Paulin}, we give a corrected statement and proof of this proposition.

In Section 4 we present \Cref{treealg}, and prove that it terminates after finitely many steps. We discuss its implementation and give some examples which compare and contrast it with the algorithm from \cite{SL2}. 

In Section 5, we show that the same method can be applied to determine whether or not two-generated subgroups of the isometry group of a locally finite simplicial tree are free and discrete, with respect to the topology of pointwise convergence (which, in this setting, is equivalent to the compact-open topology). For any of these subgroups (including those of ${\rm SL_2}(K)$) which are discrete and free, we show that there is also a practical algorithm to solve the constructive membership problem.

\section{Background}

A \textit{local field} is a field which is locally compact with respect to the topology induced by some non-trivial absolute value. Such a field $K$ is said to be \textit{non-archimedean} if the corresponding absolute value $|-|$ is non-archimedean, meaning it satisfies the \textit{ultrametric inequality}
\begin{align*}
|a+b|\le \max\{|a|,|b|\},
\end{align*}
for all $a, b \in K$. We note that equality holds when $|a|\ne |b|$. 

Any local field that does not satisfy the ultrametric inequality is said to be \textit{archimedean}, and is isomorphic to either $\rr$ or $\cc$ with the same topology as that induced by the standard absolute values; see \cite[Chapter 3, Theorem 1.1]{LF}. Non-archimedean local fields are a little different, and have an equivalent characterisation in terms of valuations. 

A \textit{valuation} on a field $K$ is a group homomorphism $v\colon K^\times \to \rr$ such that, when extended by defining $v(0)=\infty$, the ultrametric inequality holds for all $x,y \in K$:
\begin{align*}
v(x+y)\ge \min\{v(x),v(y)\}.
\end{align*}
We say that $v$ is \textit{discrete} if $v(K^\times)\cong \zz$. Given any valuation $v$ on a field $K$, the \textit{ring of integers} $\oo=\{x \in K : v(x) \ge 0\}$ is a principal ideal domain with unique maximal ideal $\pp=\{x \in K : v(x) > 0\}$. The quotient $k=\oo/\pp$ is called the \textit{residue field} of $K$. Furthermore, setting $|x|_v=c^{-v(x)}$ for some $c \in (1,\infty)$ defines a non-archimedean absolute value on $K$. A field $K$, equipped with discrete valuation $v$, that is complete with respect to $|-|_v$ and has finite residue field $k$ is a non-archimedean local field. The converse also holds, giving two equivalent definitions of a non-archimedean local field; see \cite[Chapter 4]{LF} for further details.

For a non-archimedean local field $K$, the maximal ideal $\pp$ is generated by a \textit{uniformiser} $\pi \in \oo$ such that $v(\pi)=1$, and hence the residue field $k$ is of the form $\oo/\pi \oo$. For a fixed finite set $S$ of representatives of $k$, every $a\in K^\times$ can be uniquely expressed a sum
\begin{align*}
a=\sum\limits_{i=N}^\infty a_i\pi^i, 
\end{align*}
with each $a_i \in S$, and for some integer $N$ such that $a_N \ne 0$; see \cite[Chapter 4]{LF}. It follows that non-archimedean local fields satisfy the \textit{Bolzano-Weierstrass property}, that is, every bounded sequence (in terms of the corresponding absolute value) has a convergent subsequence.

A common example of a non-archimedean local field is the $p$-adic numbers, defined using the \textit{$p$-adic valuation} $v_p$ on $\qq$. Namely, if $p$ is a prime and $x\in \qq$ is of the form $p^r\frac{a}{b}$ with $p \nmid a,b$, then $v_p(x)=r$. The corresponding absolute value is usually defined by $|x|_p=p^{-r}$, and the \textit{$p$-adic numbers} $\qq_p$ are the completion of $\qq$ with respect to $|-|_p$. Every non-archimedean local field is isomorphic to a finite extension of either $\qq_p$ or the field of formal Laurent series $\mathbb{F}_p((t))$ for some prime $p$; see \cite[Exercise 25 of Chapter 4 and Lemma 1.1 of Chapter 8]{LF}.

\medskip
 
Given a non-archimedean local field $K$ with associated valuation $v$, there is a locally finite simplicial tree $T_v$, called the \textit{Bruhat-Tits tree}, upon which the group ${\rm SL_2}(K)$ acts. The vertices of $T_v$ are equivalence classes of free $\oo$-modules of rank two (called \textit{lattices}), where lattices $L$ and $L'$ are \textit{equivalent} if $L=xL'$ for some $x\in K^\times$. Furthermore, given a lattice $L$, each equivalence class of lattices has a unique representative $L_0\subseteq L$ for which $L/L_0$ is isomorphic (as an $\oo$-module) to $\oo/\pi^n\oo$, for some $n\in \zz_{\ge 0}$. This gives rise to the edge structure of $T_v$, by having edges between the vertices represented by $L$ and $L_0$ if and only if $n=1$; for further details, see \cite[Chapter II]{Serre}.

There is a natural action of ${\rm GL_2}(K)$ on the set of lattices, and this gives rise to a faithful action of ${\rm PGL_2}(K)$ on $T_v$ by isometries. Moreover, the subgroups ${\rm SL_2}(K)$ and ${\rm PSL_2}(K)$ act on $T_v$ \textit{without inversions}, that is, no element swaps adjacent vertices; see \cite[Corollary II.3.14]{MS}. Isometries of a simplicial tree $T$ acting without inversions can be classified based on their \textit{translation length}: given such an isometry $g$, this is the integer
$$l(g)=\min_{x\in V(T)}d(x,gx),$$
where $V(T)$ denotes the vertex set of $T$, and $d$ is the standard path metric on $T$. Note that $l(g)=l(g^{-1})$ and $l(hgh^{-1})=l(g)$ for all such isometries $g, h$ of $T$. Moreover, if $l(g)=0$, then $g$ fixes a vertex of $T$ and $g$ is said to be \textit{elliptic}. If $l(g)>0$ then $g$ is said to be \textit{hyperbolic}.

\begin{prop}\label{tl}
Suppose that $g$ is a hyperbolic isometry of a simplicial tree $T$. Then $\{p \in V(T) : d(p,gp)=l(g)\}$ is the vertex set of a straight path in $T$ (called the \textup{axis} of $g$) on which $g$ acts by translations of length $l(g)$. Moreover, if a vertex $q\in V(T)$ is at distance $k$ from the axis of $g$, then $d(q,gq)=l(g)+2k$.
\end{prop}
\begin{proof}
See \cite[Chapter I, Proposition 24]{Serre}.
\end{proof}

\begin{cor}\label{edgeaxis}
An edge $p-q$ in $T$ is contained in the axis of a hyperbolic element $g$ if and only if $d(p,gp)=d(q,gq)$.
\end{cor}

Elements of ${\rm SL_2}(K)$ can be classified as either elliptic or hyperbolic via their action on the Bruhat-Tits tree $T_v$, and this depends only on the trace:

\begin{prop}\label{tv}
If $A \in {\rm SL_2}(K)$, then $l(A)=-2\min\{0,v({\rm tr}(A))\}$.
\end{prop}
\begin{proof}
See \cite[Proposition II.3.15]{MS}.
\end{proof}

\section{Discrete and free subgroups}

In this section we fix a non-archimedean local field $K$ with valuation $v$, and present key results which underpin our algorithm that determines whether or not a given two-generated subgroup of ${\rm SL_2}(K)$ is discrete and free of rank two. As with the algorithm for two-generated subgroups of ${\rm SL_2}(\rr)$ in \cite{SL2}, we use Nielsen transformations on pairs of generating elements, but in this case we aim to minimise translation lengths until either an elliptic element or a suitable pair of hyperbolic elements is encountered (a similar `reduction' process is used in Section 4 of \cite{CV} in the context of free groups of rank two acting on $\rr$-trees). We also show that a group containing an elliptic element cannot be both discrete and free, and give some translation length formulae which allow us to check when a pair of hyperbolic elements generate a discrete and free group of rank two.

\medskip

First recall that a \textit{Nielsen transformation} takes an $n$-tuple of elements $(g_1, \dots, g_n)$ of a group and performs some finite sequence of the following operations:
\begin{itemize}
\item Swap $g_i$ and $g_j$ (for $i \ne j$);
\item Replace $g_i$ by $g_i^{-1}$;
\item Replace $g_i$ by $g_j^{-1}g_i$ (for $i \ne j$).
\end{itemize}
This preserves generation of the subgroup generated by $g_1, \dots, g_n$.  

Recall also that a \textit{topological group} is a group equipped with a topology such that the inversion and multiplication maps are continuous. A topological group is said to be \textit{discrete} if the corresponding topology is discrete. Since multiplication by any element is a homeomorphism, such a group is discrete if and only if the set $\{1\}$ is open. Hence any metrisable topological group (in particular, ${\rm SL_2}(K)$ - via the subspace topology and metric it inherits from $K^4$) is discrete if and only if any sequence of elements in the group converging to the identity is eventually constant. 

\begin{prop}\label{disck}
Let $A \in {\rm SL_2}(K)$. Then the subgroup $\langle A \rangle \le {\rm SL_2}(K)$ is discrete if and only if either $A$ has finite order or $v({\rm tr}(A))<0$.
\end{prop}
\begin{proof}
Set $A=\left[ \begin{array}{cc}
a & b \\
c & d \end{array} \right]$ and $t={\rm tr}(A)$. If $A$ has finite order then it generates a discrete group, so suppose that $v(t)<0$, that is, $|t|_v>1$. Using the ultrametric inequality, we may also assume that $|a|_v>1$. Let $a_n$ denote the top left entry of the matrix $A^n$ for each $n \in \nn$. By the Cayley-Hamilton Theorem we have $A^n=tA^{n-1}-A^{n-2}$ so, if $|a_{n-1}t|_v>|a_{n-2}|_v$, then the ultrametric inequality implies
\begin{align*}
|a_nt|_v>|a_{n-1}t-a_{n-2}|_v=|a_{n-1}t|_v>|a_{n-1}|_v.
\end{align*}
Since $|a_1t|_v>1=|a_0|_v$, this inductively proves that $|a_nt|_v>|a_{n-1}|_v$ and hence that $|a_{n+1}|_v=|a_nt|_v$ for all $n \in \nn$. Thus $|a_n|_v$ tends to $\infty$ as $n$ does, so $\langle A \rangle$ is discrete. 

On the other hand, suppose $A$ has infinite order and $v(t)\ge 0$, that is, $|t|_v\le 1$. Let $a_n, b_n, c_n$ and $d_n$ denote the corresponding entries of the matrix $A^n$. Note that if both $|a_{n-1}|_v$ and $|a_{n-2}|_v$ are bounded above, then so is $|a_n|_v$ by the ultrametric inequality and the Cayley-Hamilton Theorem. It follows by induction that $|a_n|_v$ is bounded above for all $n \in \nn$. Similarly, $|b_n|_v, |c_n|_v$ and $|d_n|_v$ are bounded above for all $n \in \nn$. The Bolzano-Weierstrass property then implies that $\langle A \rangle$ is not discrete.
\end{proof}

\begin{cor}\label{noell}
If $G \le {\rm SL_2}(K)$ is discrete and free then $l(g)>0$ for all $g \in G$.
\end{cor}
\begin{proof}
Suppose that $g \in G$ is elliptic. Then either $g$ has finite order, whereby $G$ is not free, or otherwise \Cref{tv} implies that $v({\rm tr}(A))\ge 0$. But then $G$ cannot be discrete by \Cref{disck}.
\end{proof}

We will frequently make use of the following version of the Ping Pong Lemma. As stated, it applies only to metrisable topological groups acting continuously on a topological space; this makes it more specialised than other statements of the lemma (for instance, see \cite{SL2} or \cite{LU}) but it enables us to determine when such a group is not only free, but discrete as well. 

Recall that a topological group $G$ acts \textit{continuously} on a topological space $X$ if the map $G\times X \to X$ (given by $(g,x) \mapsto gx$) is continuous with respect to the product topology. Note that the action of ${\rm SL_2}(K)$ on the Bruhat-Tits tree $T_v$ is defined by polynomials and is hence continuous. 

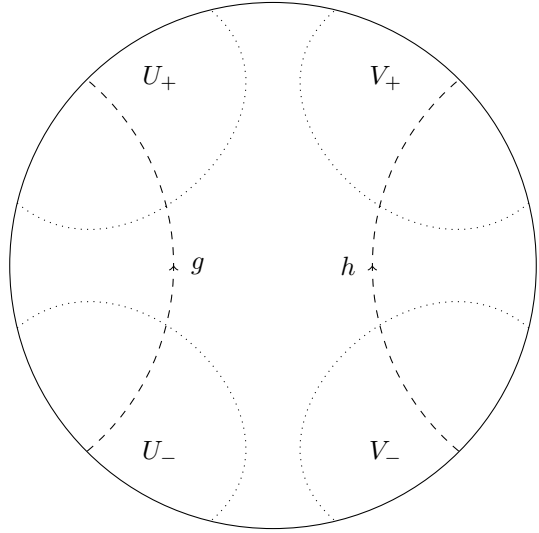
\begin{figure}[h]
\centering
\begin{tikzpicture}
  [scale=1,auto=left] 

\draw (0,0) circle (3.5cm);

\draw [dotted] (3.4,0.83) to [out=220,in=220, distance=2 cm] (0.83,3.4) ;
\draw [dotted] (3.4,-0.83) to [out=140,in=140, distance=2 cm] (0.83,-3.4) ;
\draw [dotted] (-3.4,0.83) to [out=320,in=320, distance=2 cm] (-0.83,3.4) ;
\draw [dotted] (-3.4,-0.83) to [out=40,in=40, distance=2 cm] (-0.83,-3.4) ;

\draw [dashed, ->-] (2.47,-2.47) to [out=140,in=220, distance=2 cm] (2.47, 2.47);
\draw [dashed, ->-] (-2.47,-2.47) to [out=40,in=320, distance=2 cm] (-2.47, 2.47);

 \node at (-1,0) {$g$};
  \node at (1,0) {$h$};

 \node at (-1.5,2.5) {$U_+$};
  \node at (-1.5,-2.5) {$U_-$};
   \node at (1.5,2.5) {$V_+$};
  \node at (1.5,-2.5) {$V_-$};

\end{tikzpicture}
\caption{The Ping Pong Lemma} \label{PPL}
\end{figure}

\begin{lem}[The Ping Pong Lemma]\label{PPLemma}
Let $G$ be a metrisable topological group acting continuously on a topological space $X$ and let $g,h \in G\backslash\{1\}$. Suppose that $U_+, U_-, V_+, V_-$ are non-empty closed pairwise disjoint subsets of $X$ which do not cover $X$ and satisfy 
\begin{align*}
& g(X\backslash U_-) \subseteq U_+; &g^{-1}(X\backslash U_+) \subseteq U_-; \\
& h(X\backslash V_-) \subseteq V_+; &h^{-1}(X\backslash V_+) \subseteq V_-.
\end{align*}
Then the subgroup $H=\langle g,h \rangle\le G$ is discrete and free of rank two.
\end{lem}
\begin{proof}
Fix some $x\in D=X\backslash (U_+\cup U_-\cup V_+\cup V_-)\neq \varnothing$. If $w\in H$ is a non-trivial word in $g,h$ then note that $w(x)\in X\backslash D$ by hypothesis. In particular, this implies $w\neq 1$ in $H$ and thus $H$ is free of rank two. On the other hand, suppose that $H$ is not discrete. Then one can find a sequence $(h_n)_{n\in\nn}$ of non-identity elements of $H$ which converges to $1 \in H$. Since $h_n(x) \in X\backslash D$ for each $n \in \nn$ and $G$ acts continuously on $X$, this gives a sequence $(h_n(x))_{n\in\nn}$ of elements of $X\backslash D$ which converges to $x\in D$. But $X\backslash D$ is closed, so this is impossible. Thus $H$ is discrete and free of rank two. See \Cref{PPL}.
\end{proof}

\begin{figure}[h]
\centering
\begin{tikzpicture}
  [scale=1,auto=left] 
     \draw [dashed] (-4,-1.5) to (-4,-0.5) ;  \draw [dashed] (-4,1.5) to (-4,-0.5) ;
                  \node[circle,inner sep=0pt,minimum size=3,fill=black] (1) at (-4,1.5) {};
                  \node[circle,inner sep=0pt,minimum size=3,fill=black] (1) at (-4,0) {};
                       \draw [dashed,->] (-4,1.5) to (-4,3);
                            \node[circle,inner sep=0pt,minimum size=3,fill=black] (1) at (-4,-1.5) {};
                   \draw [dashed,>-]  (-4,-3) to (-4,-1);
      \node at (-4.2,-1.5) {$p$};  \node at (-4.2,0) {$p'$};  \node at (-4.3,1.65) {$Ap$};
       \node at (-4, -3.3) {${\rm Axis}(A)$};  \node at (-4, 3.3) {${\rm Axis}(A)$};
             
            \draw [dashed] (-1.5,-1.5) to (-1.5,-0.5) ; \draw [dashed] (-1.5,-0.5) to (-1.5,1.5) ;
                           \node[circle,inner sep=0pt,minimum size=3,fill=black] (1) at (-1.5,1.5) {};
                            \node[circle,inner sep=0pt,minimum size=3,fill=black] (1) at (-1.5,0) {};
                                         \draw [dashed,->]  (-1.5,1.5) to (-1.5,3);
                            \node[circle,inner sep=0pt,minimum size=3,fill=black] (1) at (-1.5,-1.5) {};
                             \draw [dashed,>-] (-1.5,-3) to (-1.5,-1.5);
     
   \node at (-1.25,-1.5) {$q$};  \node at (-1.25,0) {$q'$};\node at (-1.15,1.65) {$Bq$};
            \node at (-1.5, -3.3) {${\rm Axis}(B)$};     \node at (-1.5, 3.3) {${\rm Axis}(B)$};

  \draw [dotted] (-5,2.5) to [out=285,in=255, distance=2 cm] (-3,2.5) ;
       \node at (-3.6,2.2) {$U_+$};
         \draw [dotted] (-2.5,2.5) to [out=285,in=255, distance=2 cm] (-0.5,2.5) ;
       \node at (-1.9,2.2) {$V_+$};     
         \draw [dotted] (-5,-2.5) to [out=70,in=110, distance=2 cm] (-3,-2.5) ;
       \node at (-3.6,-2.2) {$U_-$};    
                \draw [dotted] (-2.5,-2.5) to [out=70,in=110, distance=2 cm] (-0.5,-2.5) ;
       \node at (-1.9,-2.2) {$V_-$};

\draw [dashed, |-|] (3.3,1) to (3.3,-1); \node at (4,0) {$\Delta(A,B)$};

 \node[circle,inner sep=0pt,minimum size=3,fill=black] (1) at (3,-1) {};
     \draw [dashed] (3,-1) to (3,1); 
                  \node[circle,inner sep=0pt,minimum size=3,fill=black] (1) at (3,1) {};
                       \draw [dashed]  (3,1) to (2,1.5) ;
                  \node[circle,inner sep=0pt,minimum size=3,fill=black] (1) at (2,1.5) {};
                                         \draw [dashed]  (3,1) to (4,1.5) ;
                  \node[circle,inner sep=0pt,minimum size=3,fill=black] (1) at (4,1.5) {};
                                \draw [dashed,->]  (4,1.5) to (4.25,3);
                                  \draw [dashed,->]  (2,1.5) to (1.75,3);
 \draw [dashed] (3,-1) to (2,-1.5) ;
                  \node[circle,inner sep=0pt,minimum size=3,fill=black] (1) at (4,-1.5) {};
                                         \draw [dashed] (3,-1) to (4,-1.5) ;
                  \node[circle,inner sep=0pt,minimum size=3,fill=black] (1) at (2,-1.5) {};
                                \draw [dashed,>-]   (1.75 ,-3) to (2,-1.5);
                                  \draw [dashed,>-]   (4.25,-3) to (4,-1.5);
                                  
       \node at (1.8,-1.5) {$p$};   \node at (1.7,1.65) {$Ap$};                               
        \node at (4.2,-1.5) {$q$};   \node at (4.3,1.65) {$Bq$};                         
                                  
                                        \node at (1.75, 3.3) {${\rm Axis}(A)$}; \node at (1.75, -3.3) {${\rm Axis}(A)$};
                                         \node at (4.25, 3.3) {${\rm Axis}(B)$};  \node at (4.25, -3.3) {${\rm Axis}(B)$};

  \draw [dotted] (1.1,2.5) to [out=285,in=255, distance=2 cm] (2.9,2.5) ;
       \node at (2.3,2.2) {$U_+$};
       
         \draw [dotted] (3.1,2.5) to [out=285,in=255, distance=2 cm] (4.9,2.5) ;
       \node at (3.8,2.2) {$V_+$};
       
                \draw [dotted] (1.1,-2.5) to [out=70,in=110, distance=2 cm] (2.9,-2.5) ;
       \node at (2.3,-2.2) {$U_-$};

                \draw [dotted] (3.1,-2.5) to [out=70,in=110, distance=2 cm] (4.9,-2.5) ;
       \node at (3.8,-2.2) {$V_-$};

\end{tikzpicture} 
\caption{Applying the Ping Pong Lemma on trees} \label{TreePPL}
\end{figure}
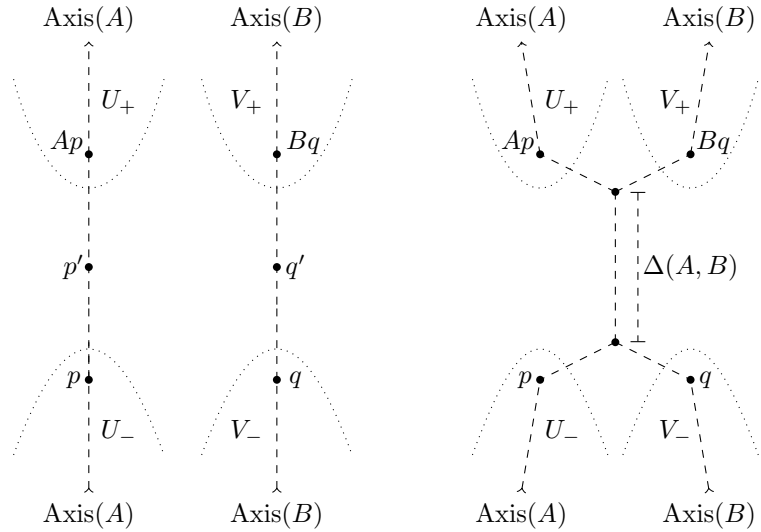

Using a version of the Ping Pong Lemma that does not involve discreteness, Lemma 2.6 of \cite{CM} shows that two hyperbolic isometries of a $\rr$-tree generate a free group of rank two when their axis overlap is sufficiently small. Lemma 3.2 of \cite{UZ} generalises this to $\Lambda$-trees (where distances take values in some totally ordered abelian group $\Lambda$, not necessarily $\rr$ or $\zz$). Here we use our version of the Ping Pong Lemma to prove a similar, but stronger, result for certain hyperbolic isometries of a simplicial tree:

\begin{prop}\label{hyptree}
Let $G$ be a metrisable topological group acting continuously, by isometries and without inversions on a simplicial tree $T$. Suppose that $A, B \in G$ are hyperbolic, and their axes are either disjoint or intersect along a path of length $0\le\Delta(A,B)<\min\{l(A),l(B)\}$. Then the subgroup $\langle A, B \rangle\le G$ is discrete and free of rank two.
\end{prop}
\begin{proof}
First of all, if the axes of $A$ and $B$ are disjoint, then there is a unique path $P$ of minimal distance from a vertex $p'$ on the axis of $A$ to a vertex $q'$ on the axis of $B$. Choose vertices $p$ and $q$ (on the axes of $A$ and $B$ respectively) so that the interior of the path between $p$ and $Ap$ contains $p'$, and the interior of the path between $q$ and $Bq$ contains $q'$ (if either $A$ or $B$ has translation length one, then it may be necessary to subdivide each edge of $T$ at its midpoint in order to find such vertices); see the left-hand diagram of \Cref{TreePPL}. On the other hand, if the axes of $A$ and $B$ intersect along a common subpath $P$ of length $\Delta(A,B) < \min\{l(A),l(B)\}$, then choose vertices $p$ and $q$ (on the axes of $A$ and $B$ respectively) such that the interior of the paths between $p$ and $Ap$, and $q$ and $Bq$, each contain $P$ (if either $A$ or $B$ has translation length $\Delta(A,B)+1$, then it may be necessary to subdivide each edge of $T$ at its midpoint in order to find such vertices); see the right-hand diagram of \Cref{TreePPL}.

In each case, define $U_+$ (respectively $U_-$) to be the maximal subtree of $T$ containing all vertices on the axis of $A$ from $Ap$ onwards (respectively up to, and including $p$) with respect to the direction of translation, but no other vertices on the axis of $A$. Similarly define $V_+$ (respectively $V_-$) as the maximal subtree containing the vertices of the axis of $B$ from $Bq$ onwards (respectively up to, and including, $q$) but no other vertices on the axis of $B$. Then, in each case, $U_-,U_+, V_-$ and $V_+$ are non-empty, pairwise disjoint closed subsets that do not cover $T$. Moreover, \Cref{tl} implies that $A(T\backslash U_-) \subseteq U_+$, $A^{-1}(T\backslash U_+) \subseteq U_-$, $B(T\backslash V_-) \subseteq V_+$ and $B^{-1}(T\backslash V_+) \subseteq V_-$. The result then follows from the Ping Pong Lemma.
\end{proof}

Given two hyperbolic isometries $A$ and $B$ of a simplicial tree, determining how their axes interact relies on the following proposition. It is effectively a reformulation of \cite[Proposition 1.6]{Paulin} for isometries of simplicial trees, however we provide an extra case (given by case $(2)(iii)$ in our version of the proposition) which was not considered in \cite{Paulin}. In the appendix, we give an erratum to \cite[Proposition 1.6]{Paulin} with the author of \cite{Paulin} and prove this extra case in the context of $\rr$-trees.

\begin{prop}\label{overlap}
Let $A$ and $B$ be hyperbolic isometries of a simplicial tree, such that $AB$ and $A^{-1}B$ act without inversions. Then precisely one of the following holds:
\begin{enumerate}[label={$(\arabic*)$}]
\item The axes of $A$ and $B$ do not intersect. If $k$ is the minimum distance between the two axes, then 
$$l(AB)=l(A^{-1}B)=l(A)+l(B)+2k.$$
\item The axes of $A$ and $B$ intersect along a (possibly infinite) path of length $\Delta=\Delta(A,B)\ge 0$, 
$$\max\{l(AB), l(A^{-1}B)\}=l(A)+l(B),$$
and either:
\begin{enumerate}[label={$(\roman*)$}]
\item $\Delta<\min\{l(A),l(B)\}$ and $\min\{l(AB), l(A^{-1}B)\}=l(A)+l(B)-2\Delta$; \\ or
\item $\Delta>\min\{l(A),l(B)\}$ and $\min\{l(AB), l(A^{-1}B)\}=|l(A)-l(B)|$; \\ or
\item $\Delta=\min\{l(A),l(B)\}$, either the axes of $B$ and $A^{-1}BA$ (if $l(A) \le l(B)$) or the axes of $A$ and $B^{-1}AB$ (if $l(A) > l(B)$) intersect along a (possibly infinite) path of length $\Delta' \ge 0$ and \\
$$\hspace{2cm} \min\{l(AB), l(A^{-1}B)\}=\left\{
\begin{array}{ll}
|l(A)-l(B)|-2\Delta' & \textup{ if } \Delta' < \frac{|l(A)-l(B)|}{2} \\
0 & \textup{ otherwise. }\\
\end{array} 
\right.$$
\end{enumerate}
\end{enumerate}
\end{prop}
\begin{proof}
\Cref{overlap} follows from \Cref{prop:appendix}, with essentially the same proof. The only difference is that the proof of the third subcase of \Cref{prop:appendix} $(2)(ii)$ (corresponding to the second case of \Cref{overlap} $(2)(iii)$) uses the fact that an isometry of an $\rr$-tree which fixes the midpoint $m$ of some path is elliptic. In the context of simplicial trees, however, this midpoint $m$ could be a vertex or the midpoint of an edge. One can check that the assumption that both $AB$ and $A^{-1}B$ act without inversions is sufficient to ensure that this midpoint $m$ is indeed a vertex and thus $\min\{l(AB), l(A^{-1}B)\}=0$, as desired.
\end{proof}

We note that the missing case from \cite[Proposition 1.6]{Paulin} was discovered when considering various examples in ${\rm SL_2}(\qq_7)$. Namely, given the matrices
$$X=\left[ \begin{array}{cc}
7^3 & 0 \\ [6pt]
 0 & \frac{1}{7^3}  \end{array} \right], 
Y=\left[ \begin{array}{cc}
\frac{2}{7^7} & 7^3 \\ [6pt]
\frac{1}{7^3} &  7^7 \end{array} \right],$$
setting $A=XY$ and $B=X^3Y^3$ yields hyperbolic elements with respective translation lengths of 8 and 32. Moreover, the axes of $A^{-1}$ and $B$ overlap with opposite directions of translation. However $l(A^{-1}B)=16$, and this is inconsistent with the formula given in case $(2)(ii)$ of \cite[Proposition 1.6]{Paulin}; this value is neither $l(B)-l(A)$ nor of the form $l(A)+l(B)-2\Delta$ for some $\Delta<8$.

\begin{cor}\label{compare}
Let $G$ be a metrisable topological group acting continuously, by isometries and without inversions on a simplicial tree. If $A, B\in G$ are hyperbolic and $|l(A)-l(B)|<\min\{l(AB), l(A^{-1}B)\}$, then $\langle A, B \rangle\le G$ is discrete and free of rank two.
\end{cor}
\begin{proof}
We consider the cases given in \Cref{overlap}. If the axes of $A$ and $B$ do not intersect then $$l(AB)=l(A^{-1}B)\ge l(A)+l(B)>|l(A)-l(B)|.$$
If the axes of $A$ and $B$ do intersect, and $\Delta(A,B) < \min\{l(A),l(B)\}$, then $$\min\{l(AB), l(A^{-1}B)\}=l(A)+l(B)-2\Delta(A,B) > |l(A)-l(B)|.$$
Otherwise, we have $\Delta(A,B) \ge \min\{l(A),l(B)\}$ and $$\min\{l(AB), l(A^{-1}B)\} \le |l(A)-l(B)|.$$
Hence $|l(A)-l(B)|<\min\{l(AB), l(A^{-1}B)\}$ if and only if the axes of $A$ and $B$ either do not intersect, or intersect along a path of length $0\le \Delta(A,B) < \min\{l(A),l(B)\}$. By \Cref{hyptree}, this implies $\langle A, B \rangle \le G$ is discrete and free of rank two.
\end{proof}

We conclude this section by noting that determining whether or not a finitely generated subgroup of ${\rm SL_2}(K)$ is discrete and free is equivalent to the same problem for the corresponding subgroup in ${\rm PSL_2}(K)$ (which inherits the quotient topology from ${\rm SL_2}(K)$).

\begin{prop}\label{psl}
Let $K$ be a local field and suppose $G \le {\rm SL_2}(K)$ is $n$-generated. Then $G$ is discrete and free of rank $n$ if and only if the corresponding subgroup $\overline{G}\le {\rm PSL_2}(K)$ (its image under the quotient map) is discrete and free of rank $n$.
\end{prop}
\begin{proof}
It is easy to check that $G$ is discrete if and only if $\overline{G}$ is. So consider the restriction of the quotient map $\pi\colon {\rm SL_2}(K) \to {\rm PSL_2}(K)$ to the epimorphism $\pi_G \colon G \to \overline{G}$. Note that $\pi(g)=1$ if and only if $g=\pm I_2$. So if $G$ is free of rank $n$ then $\pi_G$ is 1-to-1. Thus $G \cong \overline{G}$ and so $\overline{G}$ must also be free of rank $n$. 

Similarly, if $\overline{G}$ is free of rank $n$ then, by the universal property of free groups, there exists a unique homomorphism $\overline{G} \to G$ sending the generators of $\overline{G}$ back to their corresponding elements in $G$. This is an inverse to $\pi_G$, showing that $G \cong \overline{G}$ and so $G$ must also be free of rank $n$.
\end{proof}

\section{The algorithm}

In this section we present our algorithm, which determines after finitely many steps whether or not a two-generated subgroup of ${\rm SL_2}(K)$ is discrete and free of rank two. The key idea is to use \Cref{tv} to compute translation lengths on the Bruhat-Tits tree, and perform Nielsen transformations on the generators until these produce either an elliptic element, or two hyperbolic elements satisfying the hypotheses of \Cref{compare}. By \Cref{psl}, the algorithm can also be applied to two-generated subgroups of ${\rm PSL_2}(K)$ by taking representatives in ${\rm SL_2}(K)$.

\begin{alg}\label{treealg}
Let $K$ be a non-archimedean local field. Given two elements $A, B \in {\rm SL_2}(K)$, we proceed as follows. If $G=\langle A, B \rangle \le {\rm SL_2}(K)$ is discrete and free of rank two then the algorithm will return true and output a generating pair for $G$ which satisfy the hypotheses of the Ping Pong Lemma; otherwise it will return false. 
\begin{enumerate}[label={$(\arabic*)$}]
\item Set $X=A$, $Y=B$. If $l(X)=0$ or $l(Y)=0$ then return false.
\item If $l(X)> l(Y)$ then swap $X$ and $Y$.
\item Compute $m=\min\{l(XY), l(X^{-1}Y)\}$.
\item If $m=0$ then return false.
\item If $m\le l(Y)-l(X)$ then replace $Y$ by an element from $\{XY, X^{-1}Y\}$ which has translation length $m$ and return to $(2)$.
\item Otherwise return true and the generating pair $(X,Y)$.
\end{enumerate}
\end{alg}

\begin{thm}\label{pf1}
\Cref{treealg} terminates after finitely many steps and produces the correct output.
\end{thm}
\begin{proof}
If at any point the algorithm encounters an elliptic element then $G$ is not discrete and free by \Cref{noell}. So suppose that the algorithm only ever encounters hyperbolic elements. Then it must reach step $(5)$. If $m>l(Y)-l(X)$ then, by \Cref{compare}, $G$ is discrete and free and the elements $X$ and $Y$ satisfy the hypotheses of the Ping Pong Lemma. Otherwise the algorithm performs a Nielsen transformation, and outputs a new pair of generators for $G$ on which to run the algorithm.

If this sequence of Nielsen transformations never terminates, then there is an infinite sequence $(x_n,y_n)=(l(X_n),l(Y_n))$ of integral translation length pairs which satisfies $0<x_n\le y_n$ for all $n\in\nn$ and is decreasing in each component; such a sequence must converge. Moreover, for each pair $(X_n, Y_n)$ of generators, we are in either case $(2)(ii)$ or the first subcase of $(2)(iii)$ of \Cref{overlap}. Hence, at each stage $(x_n, y_n)$ is replaced by either $(y_n-x_n-k_n,x_n)$ or $(x_n, y_n-x_n-k_n)$ for some $0\le k_n < y_n-x_n$. In particular, this implies that $x_{n+1}+y_{n+1}=y_n-k_n$ for all $n\in \nn$. Rearranging and taking limits, it follows that $\lim\limits_{n\to\infty}x_n=-\lim\limits_{n\to\infty}k_n\le 0$, a contradiction since each $x_n$ is a positive integer. Hence this algorithm must eventually terminate, proving the theorem.
\end{proof}
In terms of implementing this algorithm in a computational package such as {\sc magma}, the software needs to be able to perform matrix multiplications over $K$, and compute traces and valuations. Since each non-zero element of $K$ can be expressed uniquely in the form $\sum\limits_{i=N}^\infty a_i\pi^i$ for some integer $N$ with $a_N \neq 0$ and some uniformiser $\pi$, computing valuations and performing both addition and multiplication over $K$ is straightforward. But there is a clear obstacle in the computational storage space needed for elements of $K$ with an infinite expression of the above form. This can theoretically be overcome by storing elements of $K$ in terms of the data $\lbrace\pi; a_N, a_{N+1}, \dots, a_M\rbrace$ up to some appropriate finite $M$. 

Indeed, given matrices $$A=\left[ \begin{array}{cc}
a & b \\
c & d \end{array} \right], 
B=\left[ \begin{array}{cc}
e & f \\
g & h \end{array} \right],$$
one iteration of \Cref{treealg} requires firstly computing $l(A)=-2\min\{0,v(a+d)\}$ and $l(B)=-2\min\{0,v(e+h)\}$. Since any non-negative valuation gives a translation length of 0, calculating these accurately requires storing the entries of $A$ and $B$ only up to the coefficient of $\pi^0$ (that is, $M=0$ will suffice). On the other hand, assuming that $0<l(A)\le l(B)$, the first iteration of \Cref{treealg} will also require computing $l(AB)=-2\min\{0,v(ae+bg+cf+dh)\}$ and  $l(A^{-1}B)=-2\min\{0,v(de-bg-cf+ah)\}$. Storing the entries of $A$ and $B$ up to the coefficient of $\pi^{-\min\{0, v(a),v(b), \dots, v(h)\}}$ is sufficient to compute these valuations accurately. It follows inductively that storing the $\pi^i$-coefficients of entries of $A$ and $B$ up to $M=-r\min\{0,v(a),v(b), \dots, v(h)\}$ is enough to correctly apply $r$ iterations of \Cref{treealg}. Thus, given any two elements of ${\rm SL_2}(K)$, choosing large enough $M$ (compared with $-\min\{0,v(a),v(b), \dots, v(h)\}$) allows the algorithm to run correctly; if, however, at any point the number of iterations exceeds $\frac{M}{-\min\{0,v(a),v(b), \dots, v(h)\}}$, then a higher bound $M$ will need to be chosen and the algorithm restarted.

The examples we discuss below avoid this issue entirely for the case where $K=\qq_p$ for some prime $p$. By restricting our interest to pairs of matrices in ${\rm SL_2}(\qq)$, we can perform matrix multiplication and compute traces in the usual sense, and then consider $p$-adic valuations separately. In this particular case, it is interesting to view the subgroups generated as subgroups of both ${\rm SL_2}(\qq_p)$ and ${\rm SL_2}(\rr)$, and compare the properties of each. For instance, it is a well known consequence of the Ping Pong Lemma that the matrices 
$$A=\left[ \begin{array}{cc}
1 & 2 \\
0 & 1 \end{array} \right], 
B=\left[ \begin{array}{cc}
1 & 0 \\
2 & 1 \end{array} \right]$$
generate a discrete and free subgroup of ${\rm SL_2}(\rr)$, whereas the matrices 
$$A=\left[ \begin{array}{cc}
1 & 1 \\
0 & 1 \end{array} \right], 
B=\left[ \begin{array}{cc}
1 & 0 \\
1 & 1 \end{array} \right]$$ do not. However, neither of these pairs of matrices generate a discrete and free subgroup of ${\rm SL_2}(\qq_p)$ for any prime $p$ since a matrix of the form $\left[ \begin{array}{cc}
1 & \alpha \\
0 & 1 \end{array} \right]$ or 
$\left[ \begin{array}{cc}
1 & 0 \\
\alpha & 1 \end{array} \right]$ over a non-archimedean local field is elliptic. 

One iteration of \Cref{treealg} also shows that, for any prime $p$, the matrices 
$$A=\left[ \begin{array}{cc}
p & p-1 \\ [6pt]
\frac{-1}{p} & \frac{1}{p^2} \end{array} \right],  
B=\left[ \begin{array}{cc}
\frac{2}{p^4} & p^3 \\ [6pt]
\frac{1}{p^3} &  p^4 \end{array} \right]$$
generate a subgroup of ${\rm SL_2}(\qq_p)$ which is discrete and free of rank two. Using the same matrices as input for the algorithm in \cite{SL2} shows that they do not generate a free and discrete subgroup of ${\rm SL_2}(\rr)$ (this follows since ${\rm tr}(AB)=\frac{p+1}{p^3}<2$, so $AB$ is conjugate to a rotation matrix). On the other hand, for any prime $p\neq 2$, the matrices 
$$A=\left[ \begin{array}{cc}
p & p-1 \\ [6pt]
\frac{-1}{p} & \frac{1}{p^2} \end{array} \right],  
B=\left[ \begin{array}{cc}
\frac{2}{p^3} & p^4 \\ [6pt]
\frac{1}{p^4} &  p^3 \end{array} \right]$$
generate subgroups of both ${\rm SL_2}(\qq_p)$ and ${\rm SL_2}(\rr)$ which are discrete and free of rank two; this follows from \Cref{compare} and \cite[Theorem 4.4 (b)(iv)]{SL2} respectively.

Each of these examples requires only one iteration of \Cref{treealg}, but this is certainly not always the case. Indeed, given a prime $p \neq 2$ and positive integer $r$, the matrices
$$A=\left[ \begin{array}{cc}
p^3 & 0 \\ [6pt]
 0 & \frac{1}{p^3}  \end{array} \right], 
B=\left[ \begin{array}{cc}
\frac{2}{p^{3r+1}} & p^3 \\ [6pt]
\frac{1}{p^3} &  p^{3r+1} \end{array} \right]$$ generate a discrete and free subgroup of ${\rm SL_2}(\qq_p)$, and this requires $r+2$ iterations of \Cref{treealg}. 

\section{Generalisations and applications}

In this final section we discuss a generalisation of \Cref{treealg} to two-generator subgroups of the isometry group of any locally finite simplicial tree, and some applications to the constructive membership problem. Recall that, given a finitely generated subgroup $G=\langle g_1, \dots, g_n \rangle$ of some group $H$, and an element $h\in H$, the constructive membership problem involves determining whether or not $h$ is an element of $G$, and if it is, finding a word in $g_1, \dots, g_n$ that represents $h$.

Given any proper metric space $X$ (for instance, a locally finite simplicial tree) the isometry group ${\rm Isom}(X)$ (viewed as a subspace of $X^X$, the space of all continuous maps $X\to X$ equipped with the product topology) is a metrisable topological group; see \cite[Lemmas 5.B.3 and 5.B.5]{CH}. This topology is often known as the \textit{topology of pointwise convergence}, in the sense that a sequence $(f_i)$ in ${\rm Isom}(X)$ converges to $f\in {\rm Isom}(X)$ if and only if the sequence $(f_i(x))$ converges to $f(x)$ for each $x \in X$. Note that, for a non-archimedean local field $K$, the group ${\rm PSL_2}(K)$ (as a subgroup of the isometry group of the corresponding Bruhat-Tits tree) inherits the topology of pointwise convergence, and this coincides with the standard topology on ${\rm PSL_2}(K)$ used in this paper.

In the setting of isometry groups, the topology of pointwise convergence is equivalent to the well-known \textit{compact-open topology}; see \cite[Lemmas 5.B.1 and 5.B.2]{CH}. The pointwise convergence property of these equivalent topologies leads to an analogue of \Cref{noell} for isometries of a locally finite simplicial tree. Note that, by subdividing each edge of the tree at its midpoint, if necessary, every element of such an isometry group can be assumed to act without inversions. 

\begin{prop}\label{elltree}
Let $T$ be a locally finite simplicial tree and suppose that $G \le {\rm Isom}(T)$ is discrete (with respect to the topology of pointwise convergence) and free. Then $G$ contains no elliptic elements.
\end{prop}
\begin{proof}
Suppose $G$ contains some elliptic element $g$, which fixes some vertex $p$ of $T$. There are only finitely many vertices adjacent to $p$ and $g$ acts to permute these. This implies there is some integer $n_1$ for which $g^{n_1}$ fixes $p$ and all adjacent vertices. One continues inductively to obtain a sequence $(g^{n_i})$ of elements of ${\rm Isom}(T)$, where $g^{n_i}$ fixes all vertices at distance at most $i$ from $p$. But then $(g^{n_i}(x))$ converges to $x$ for each vertex $x$ of $T$, so $(g^{n_i})$ converges to the identity. Thus either $g$ has finite order or $G$ is not discrete.
\end{proof}

For any proper metric space $X$, the natural map ${\rm Isom}(X)\times X \to X$ is continuous; see \cite[Lemma 5.B.4 (2)]{CH}. This implies that \Cref{compare} can also be applied to the isometry group of a locally finite simplicial tree, when equipped with the topology of pointwise convergence. Thus we have the following generalisation of \Cref{treealg}:

\begin{alg}\label{isomalg}
Given two elements $A$ and $B$ in the isometry group of a locally finite simplicial tree $T$, and a method of computing translation lengths, we proceed through steps $(1)-(6)$ of \Cref{treealg}. If $G=\langle A, B \rangle \le {\rm Isom}(T)$ is discrete (with respect to the topology of pointwise convergence) and free of rank two, then the algorithm will return true and output a generating pair for $G$ which satisfies the hypotheses of the Ping Pong Lemma; otherwise it will return false. 
\end{alg}

\begin{thm}
\Cref{isomalg} terminates after finitely many steps and produces the correct output.
\end{thm}
\begin{proof}
The only difference from the proof of \Cref{pf1} is that if the algorithm encounters an elliptic element then $G$ cannot be both discrete and free by \Cref{elltree}, instead of \Cref{noell}.
\end{proof}

\Cref{isomalg} can be applied, for instance, to certain amalgamated free products. Suppose that $\Gamma=H *_C K$ is the amalgamated free product of groups $H$ and $K$ over some subgroup $C$ which is finite index in both $H$ and $K$. It is well-known that, given fixed transversals $T_H$ and $T_K$ of right coset representatives of $C$ in $H$ and $K$ respectively, each element $g \in \Gamma$ has a unique normal form
$$g=cx_1\dots x_n,$$
for some integer $n\ge 0$, where $c \in C$ and for each $i\ge 1$, either $x_i\in T_H$ and $x_{i+1} \in T_K$ or vice versa. Moreover, $\Gamma$ acts faithfully, by isometries, and without inversions, on a locally finite tree $T$ with vertices given by cosets of the form $gH$ or $gK$ and edges given by cosets $gC$, for $g \in \Gamma$; see \cite[Chapter I, Section 4]{Serre}. 

Consider the shortest normal form $cx_1\dots x_{n_0}$ of all conjugates of $g$ in $\Gamma$; such a form is \textit{cyclically reduced} in the sense that either $n_0= 0,1$ or $x_1$ and $x_{n_0}$ lie in different transversals. If $n_0$ is $0$ or $1$, then $g$ is conjugate into either $A$ or $B$ and hence $l(g)=0$. On the other hand, if $n_0>1$ then $l(g)=n_0$, which is an even integer; this follows from \cite[Lemma 2.25]{A} and \cite[Proposition 1.7]{Paulin}. Thus, given such a group $\Gamma$ and a method of computing a cyclically reduced normal form of each element (such algorithms exist since the transversals $T_H$ and $T_K$ are finite), \Cref{isomalg} can be applied to determine whether or not any two-generated subgroup of $\Gamma$ is both discrete and free.

\medskip

We conclude this paper by showing that, as is the case in \cite{SL2}, these algorithms to determine whether or not a subgroup of a certain group is both discrete and free of rank two have applications to the constructive membership problem. This requires the notion of a \textit{fundamental domain}: given a group $G$ acting on a topological space $X$, this is an open set $D\subseteq X$ such that, if $\overline{D}$ denotes the closure of $D$ in $X$, then
\begin{enumerate}[label={$(\roman*)$}]
\item $\bigcup_{g\in G}g\overline{D}=X$;
\item $gD\cap hD=\varnothing$ for all distinct $g,h \in G$.
\end{enumerate}

In the proof of \Cref{hyptree}, given a metrisable topological group $G$ (acting continuously, by isometries, and without inversions on a simplicial tree $T$) and two hyperbolic elements $A, B \in G$ whose axes are either disjoint or intersect along a sufficiently short path, we found vertices $p$ and $q$ (on the axes of $A$ and $B$ respectively) and considered their images $Ap$ and $Bq$ in order to construct subtrees $U_+, U_-, V_+, V_-\subseteq T$ satisfying the conditions of the Ping Pong Lemma; see \Cref{TreePPL}. Note that in each case $D_A$, which we define to be the interior of the path between $p$ and $Ap$ (this is isometric to an open interval in $\rr$ with integral endpoints, and is hence open in $T$), is a fundamental domain for the action of $\langle A \rangle$ on ${\rm Axis}(A)$. Similarly the open set $D_B$, defined to be the interior of the path between $q$ and $Bq$, is a fundamental domain for the action of $\langle B \rangle$ on ${\rm Axis}(B)$.

If the axes of $A$ and $B$ do not intersect, then set $D$ to be the union of $D_A$ and $D_B$ with the path between $p'$ and $q'$; otherwise, set $D=D_A\cup D_B$. Then the union of images of $\overline{D}$ under the action of $\langle A, B \rangle$ forms a subtree $S\subseteq T$ for which $D$ is a fundamental domain for the action of $\langle A, B \rangle$ on $S$; see the proof of \cite[Lemma 2.6]{CM} for further details. If one replaces the role of $T$ by this subtree $S$, then $D=T \backslash (U_+ \cup U_- \cup V_+ \cup V_-)$, where $U_+, U_-, V_+, V_-$ are as in \Cref{TreePPL}. Moreover, it follows from the proof of \Cref{hyptree} that there is at least one vertex in $D$. These observations yield the following algorithm:

\begin{alg}\label{cmp}
Given a discrete and free two-generated subgroup $G=\langle A, B \rangle$ of ${\rm SL_2}(K)$ (respectively the isometry group of a locally finite simplicial tree $T$, along with a method of computing translation lengths) and an element $C$ of the corresponding overgroup, we proceed as follows. If $C \in G$ then the algorithm will return true and output a word $w=w(a,b)$ (where $a,b$ are abstract elements generating a free group $F$ of rank two) such that $w(A,B)=C$; otherwise it will return false.
\begin{enumerate}[label={$(\arabic*)$}]
\item Run \Cref{treealg} (respectively \Cref{isomalg}) on $G$ to obtain generators $X=X(A,B), Y=Y(A,B)$ which satisfy the hypotheses of \Cref{hyptree}.
\item Replacing $T$ by an appropriate subtree, if necessary, find a fundamental domain $D=T \backslash (U_+ \cup U_- \cup V_+ \cup V_-)$ for the action of $G$ on $T$, and choose a vertex $z' \in D$.
\item Set $w=1 \in F$ and $z=Cz'$.
\item While $z \notin D$:
\begin{enumerate}[label={$(\roman*)$}]
\item If $z\in U_\pm$, then replace $z$ by $X^{\mp 1}z$ and $w$ by $wa^{\pm 1}$;
\item If $z \in V_\pm$ then replace $z$ by $Y^{\mp 1}z$ and $w$ by $wb^{\pm 1}$.
\end{enumerate}
\item If $w(X(A,B),Y(A,B))=C$ and $z=z'$ then return true and the word $w=w(X(a,b),Y(a,b))$; otherwise return false.
\end{enumerate}
\end{alg}

\begin{thm}
\Cref{cmp} terminates after finitely many steps and produces the correct output.
\end{thm}
\begin{proof}
We already know that step $(1)$ is correct and terminates after finitely many steps, and step $(2)$ is discussed in the paragraphs preceding the statement of the algorithm. The proof that the rest of the algorithm terminates after finitely many steps, and is correct, is as in \cite[Algorithm 1]{SL2}.
\end{proof}

\Cref{cmp} is another practical algorithm which can be implemented, so long as there is a method to determine whether or not a vertex lies in the fundamental domain $D$ and, if it doesn't, which of the subtrees $U_-,U_+, V_-, V_+$ it belongs to. Note that the proof of \Cref{tl} implies that, for any hyperbolic isometry $A$ of a simplicial tree $T$ and any vertex $x$ of $T$ (for instance, for the Bruhat-Tits tree $T_v$, one could take $x$ to be the vertex representing the standard lattice $\oo^2$), the midpoint of the path between $x$ and $Ax$ lies on the axis of $A$. Similarly, one can obtain a vertex on the axis of a second hyperbolic element $B$. Thus, after translating these vertices along each axis by appropriate powers of $A$ and $B$, and comparing distances between them, one should be able to obtain a rough idea of the vertices lying on each axis and hence a method of distinguishing between vertices in $U_-,U_+, V_-, V_+$ and $D$.

\section*{Acknowledgements}

The author would like to acknowledge his PhD supervisor Dr Jack Button for his advice and support, especially in helping identify the missing case from \cite[Proposition 1.6]{Paulin}. The author is also very grateful to the reviewer for their helpful comments on a previous version of this paper. This work was supported by the Cambridge and Woolf Fisher Trusts.

\appendix




\newpage
\section{The translation length of the product of hyperbolic 
isometries of $\rr$-trees}
\label{sect:appendix}
\section*{Matthew J. Conder and Fr\'ed\'eric Paulin}


\newcommand{\bprop}{\begin{prop}}
\newcommand{\eprop}{\end{prop}}
\newcommand{\ga}{\gamma}


As noticed by the first author of this appendix in the first version
of this paper, Assertion $(ii)$ of Proposition 1.6 $(2)$ in
\cite{Paulin} is incorrect. Explicit counter-examples are given after
the proof of \Cref{overlap}. This appendix serves as an erratum of the
paper \cite{Paulin} where Proposition 1.6 $(2)(ii)$ therein should be
replaced by Assertion $(2)(ii)$ of the following Proposition
\ref{prop:appendix}. Except this replacement, the remainder of the
paper \cite{Paulin} is unchanged.

The second author of this appendix is extremely grateful to the first
one for finding the mistake and for fixing it.

\medskip
We keep the notation of \cite{Paulin} in this appendix, in order to
facilitate the checking process. In particular, if $\ga$ is an
hyperbolic isometry of $T$, then $l(\ga)$ is its translation length
and $A_\ga$ is its translation axis. Most of the statements in the
following result also follow from \cite[Propositions 8.1, 8.3]{AB}.

\bprop\label{prop:appendix}  
Let $\ga,\delta$ be two hyperbolic isometries of an $\rr$-tree $T$.

\medskip
$(1)$ Assume that $A_\ga\cap A_\delta=\emptyset$. Let $D$ be the length of the
connecting arc $S$ between $A_\ga$ and $A_\delta$. Then $S$ is
contained in the translation axis of $\ga\delta$, and the isometry
$\ga\delta$ translates $S\cap A_\delta$ towards $S\cap A_\ga$. We have
$$
l(\ga\delta)=l(\ga)+l(\delta)+2D\;.
$$

$(2)$ Assume that $A_\ga\cap A_\delta\neq\emptyset$. Let $D\in [0,
+\infty]$ be the length of the intersection $A_\ga\cap A_\delta$, with
$D=0$ if this intersection is reduced to a point, and $D=\infty$ if
this intersection is noncompact.

\begin{enumerate}[label={$(\roman*)$}]
\item 
Either if $D>0$ and the translation directions of $\ga$
and $\delta$ on $A_\ga\cap A_\delta$ coincide, or if $D=0$, then
$$
l(\ga\delta)=l(\ga)+l(\delta)\;.
$$

\item 
Assume that $D>0$ and the translation directions of $\ga$
and $\delta$ are opposite on $A_\ga\cap A_\delta$. Let $D'\in[0,
+\infty]$ be the length of the (possibly empty or infinite) segment
$A_\delta\cap\,\ga A_\delta$ (resp.~$A_\ga\cap\,\delta A_\ga$) if
$l(\delta)> l(\ga)$ (resp.~ $l(\delta)<l(\ga)$), then

\medskip
$\bullet$~ $l(\ga\delta)= l(\ga)+l(\delta)-2D$ if
  $\min\{l(\ga),l(\delta)\}>D$,

\medskip
$\bullet$~ $l(\ga\delta)= |l(\ga)-l(\delta)|$ if $\min\{l(\ga),l(\delta)\}<D<
  \max\{l(\ga), l(\delta)\}$ or \\ \hspace*{0.6cm}
$\max\{l(\ga),l(\delta)\}\leq D$,

\medskip
$\bullet$~ $l(\ga\delta)= 0$ if $\min\{l(\ga),l(\delta)\}=D<
  \max\{l(\ga),l(\delta)\}\leq D+2D'$,

\medskip
$\bullet$~ $l(\ga\delta)= \max\{l(\ga),l(\delta)\}-D-2D'$ if 
$\min\{l(\ga),l(\delta)\}=D$ and \\ \hspace*{0.6cm}
$\max\{l(\ga),l(\delta)\}> D+2D'$.

\medskip
In all four cases, we have $l(\ga\delta)<l(\ga)+l(\delta)$.
\end{enumerate}
\eprop

\begin{proof}
We may assume that $l(\ga)\leq l(\delta)$. The proofs of Assertions
$(1)$ and $(2)(i)$, as well as the first two cases of Assertion
$(2)(ii)$, are the same ones as in \cite{Paulin}, see also
\cite[Propositions 8.1, 8.3]{AB}.

Hence we assume that $l(\ga)=D<l(\delta)$. In particular $D$ is finite
and nonzero, and $A_\ga\cap A_\delta$ is a compact segment which may
be written $[x,y]$ with $y=\ga x$. We denote by $z$ the point in $T$
such that $[y,z]=\ga A_\delta\cap A_\delta$, if this segment is
compact, or the point at infinity of $T$ such that $[y,z[\;=\ga
A_\delta\cap A_\delta$ otherwise.

\begin{center}
\begin{picture}(0,0)%
\includegraphics{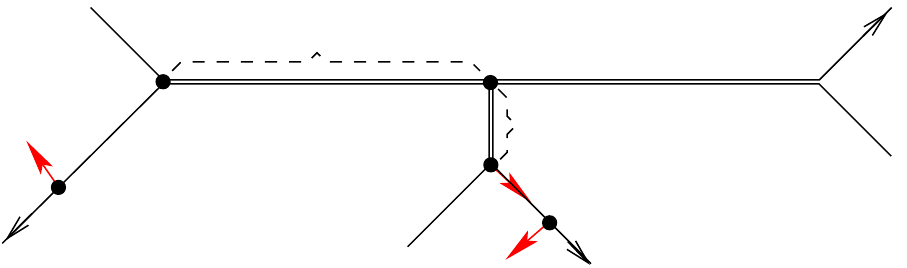}%
\end{picture}%
\setlength{\unitlength}{3812sp}%
\begingroup\makeatletter\ifx\SetFigFont\undefined%
\gdef\SetFigFont#1#2#3#4#5{%
  \reset@font\fontsize{#1}{#2pt}%
  \fontfamily{#3}\fontseries{#4}\fontshape{#5}%
  \selectfont}%
\fi\endgroup%
\begin{picture}(4503,1316)(180,-813)
\put(1621,299){\makebox(0,0)[lb]{\smash{{\SetFigFont{11}{13.2}{\rmdefault}{\mddefault}{\updefault}{\color[rgb]{0,0,0}$D=l(\ga)$}%
}}}}
\put(2772,-162){\makebox(0,0)[lb]{\smash{{\SetFigFont{11}{13.2}{\rmdefault}{\mddefault}{\updefault}{\color[rgb]{0,0,0}$D'$}%
}}}}
\put(4668,-376){\makebox(0,0)[lb]{\smash{{\SetFigFont{11}{13.2}{\rmdefault}{\mddefault}{\updefault}{\color[rgb]{0,0,0}$\ga A_\delta$}%
}}}}
\put(2454,-344){\makebox(0,0)[lb]{\smash{{\SetFigFont{11}{13.2}{\rmdefault}{\mddefault}{\updefault}{\color[rgb]{0,0,0}$z$}%
}}}}
\put(534,-475){\makebox(0,0)[lb]{\smash{{\SetFigFont{11}{13.2}{\rmdefault}{\mddefault}{\updefault}{\color[rgb]{0,0,0}$\delta z$}%
}}}}
\put(2931,-549){\makebox(0,0)[lb]{\smash{{\SetFigFont{11}{13.2}{\rmdefault}{\mddefault}{\updefault}{\color[rgb]{0,0,0}$\ga\delta z$}%
}}}}
\put(317,-708){\makebox(0,0)[lb]{\smash{{\SetFigFont{11}{13.2}{\rmdefault}{\mddefault}{\updefault}{\color[rgb]{0,0,0}$A_\delta$}%
}}}}
\put(4597,332){\makebox(0,0)[lb]{\smash{{\SetFigFont{11}{13.2}{\rmdefault}{\mddefault}{\updefault}{\color[rgb]{0,0,0}$A_\ga$}%
}}}}
\put(2634,203){\makebox(0,0)[lb]{\smash{{\SetFigFont{11}{13.2}{\rmdefault}{\mddefault}{\updefault}{\color[rgb]{0,0,0}$y=\ga x$}%
}}}}
\put(759, 61){\makebox(0,0)[lb]{\smash{{\SetFigFont{11}{13.2}{\rmdefault}{\mddefault}{\updefault}{\color[rgb]{0,0,0}$x$}%
}}}}
\end{picture}%

\end{center}

Assume first that $l(\delta)> D+2D'$, so that in particular $D'$ is
finite, $z\in T$ and $D'= d(y,z)$. See the above picture. Since
$l(\delta) > D+D'$, the point $x$ belongs to $[z,\delta z]$ and
besides $d(x,\delta z)= \ell(\delta)-D-D'> D'$. Therefore $\ga\delta
z$ does not belong to $A_\delta$. The germ at $z$ of the segment from
$z$ to $\ga\delta z$ is hence not sent to the germ at $\ga\delta z$ of
the segment from $\ga\delta z$ to $z$. Thus
$$
l(\ga\delta)= d(z,\ga\delta z)=d(\ga\delta z,y)- d(y,z)=
d(\delta z,x)- d(y,z)=\ell(\delta)-D-2D',
$$
as wanted.

\medskip
\begin{center}
\begin{picture}(0,0)%
\includegraphics{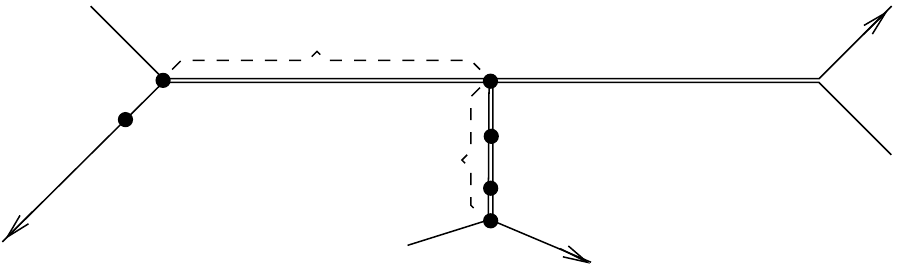}%
\end{picture}%
\setlength{\unitlength}{3812sp}%
\begingroup\makeatletter\ifx\SetFigFont\undefined%
\gdef\SetFigFont#1#2#3#4#5{%
  \reset@font\fontsize{#1}{#2pt}%
  \fontfamily{#3}\fontseries{#4}\fontshape{#5}%
  \selectfont}%
\fi\endgroup%
\begin{picture}(4503,1323)(180,-820)
\put(1621,299){\makebox(0,0)[lb]{\smash{{\SetFigFont{11}{13.2}{\rmdefault}{\mddefault}{\updefault}{\color[rgb]{0,0,0}$D=l(\ga)$}%
}}}}
\put(2161,-331){\makebox(0,0)[lb]{\smash{{\SetFigFont{11}{13.2}{\rmdefault}{\mddefault}{\updefault}{\color[rgb]{0,0,0}$D'$}%
}}}}
\put(2705,-190){\makebox(0,0)[lb]{\smash{{\SetFigFont{11}{13.2}{\rmdefault}{\mddefault}{\updefault}{\color[rgb]{0,0,0}$m$}%
}}}}
\put(4668,-376){\makebox(0,0)[lb]{\smash{{\SetFigFont{11}{13.2}{\rmdefault}{\mddefault}{\updefault}{\color[rgb]{0,0,0}$\ga A_\delta$}%
}}}}
\put(317,-708){\makebox(0,0)[lb]{\smash{{\SetFigFont{11}{13.2}{\rmdefault}{\mddefault}{\updefault}{\color[rgb]{0,0,0}$A_\delta$}%
}}}}
\put(4597,332){\makebox(0,0)[lb]{\smash{{\SetFigFont{11}{13.2}{\rmdefault}{\mddefault}{\updefault}{\color[rgb]{0,0,0}$A_\ga$}%
}}}}
\put(2634,203){\makebox(0,0)[lb]{\smash{{\SetFigFont{11}{13.2}{\rmdefault}{\mddefault}{\updefault}{\color[rgb]{0,0,0}$y=\ga x$}%
}}}}
\put(759, 61){\makebox(0,0)[lb]{\smash{{\SetFigFont{11}{13.2}{\rmdefault}{\mddefault}{\updefault}{\color[rgb]{0,0,0}$x$}%
}}}}
\put(2555,-756){\makebox(0,0)[lb]{\smash{{\SetFigFont{11}{13.2}{\rmdefault}{\mddefault}{\updefault}{\color[rgb]{0,0,0}$z$}%
}}}}
\put(838,-198){\makebox(0,0)[lb]{\smash{{\SetFigFont{11}{13.2}{\rmdefault}{\mddefault}{\updefault}{\color[rgb]{0,0,0}$\delta m$}%
}}}}
\put(2701,-466){\makebox(0,0)[lb]{\smash{{\SetFigFont{11}{13.2}{\rmdefault}{\mddefault}{\updefault}{\color[rgb]{0,0,0}$\delta^{-1} x$}%
}}}}
\end{picture}%

\end{center}

Assume now that $l(\delta)\leq D+D'$. See the above picture. Note that
$\delta^{-1} x$ does not belong to $A_\ga$ since $l(\delta)> D$, and
that $d(\delta^{-1} x,y)=l(\delta)- D\leq D'$. Let $m$ be the midpoint
of the segment $[y,\delta^{-1} x]$, so that $d(\delta m,x)= d(m,
\delta^{-1}x) =d(m,y)$.  Hence $\ga\delta m$, which is the point of
$[y,z]$ (or $[y,z[$ if $D'=+\infty$) at distance $d(\delta m,x)$ from
    $y$, is equal to $m$ and $l(\ga\delta)=0$, as wanted.

\medskip
\begin{center}
\begin{picture}(0,0)%
\includegraphics{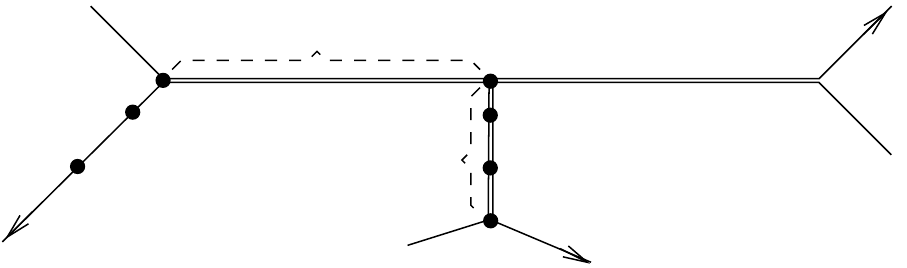}%
\end{picture}%
\setlength{\unitlength}{3812sp}%
\begingroup\makeatletter\ifx\SetFigFont\undefined%
\gdef\SetFigFont#1#2#3#4#5{%
  \reset@font\fontsize{#1}{#2pt}%
  \fontfamily{#3}\fontseries{#4}\fontshape{#5}%
  \selectfont}%
\fi\endgroup%
\begin{picture}(4503,1323)(180,-820)
\put(1621,299){\makebox(0,0)[lb]{\smash{{\SetFigFont{11}{13.2}{\rmdefault}{\mddefault}{\updefault}{\color[rgb]{0,0,0}$D=l(\ga)$}%
}}}}
\put(2161,-331){\makebox(0,0)[lb]{\smash{{\SetFigFont{11}{13.2}{\rmdefault}{\mddefault}{\updefault}{\color[rgb]{0,0,0}$D'$}%
}}}}
\put(2683,-373){\makebox(0,0)[lb]{\smash{{\SetFigFont{11}{13.2}{\rmdefault}{\mddefault}{\updefault}{\color[rgb]{0,0,0}$m$}%
}}}}
\put(4668,-376){\makebox(0,0)[lb]{\smash{{\SetFigFont{11}{13.2}{\rmdefault}{\mddefault}{\updefault}{\color[rgb]{0,0,0}$\ga A_\delta$}%
}}}}
\put(317,-708){\makebox(0,0)[lb]{\smash{{\SetFigFont{11}{13.2}{\rmdefault}{\mddefault}{\updefault}{\color[rgb]{0,0,0}$A_\delta$}%
}}}}
\put(4597,332){\makebox(0,0)[lb]{\smash{{\SetFigFont{11}{13.2}{\rmdefault}{\mddefault}{\updefault}{\color[rgb]{0,0,0}$A_\ga$}%
}}}}
\put(2634,203){\makebox(0,0)[lb]{\smash{{\SetFigFont{11}{13.2}{\rmdefault}{\mddefault}{\updefault}{\color[rgb]{0,0,0}$y=\ga x$}%
}}}}
\put(2555,-756){\makebox(0,0)[lb]{\smash{{\SetFigFont{11}{13.2}{\rmdefault}{\mddefault}{\updefault}{\color[rgb]{0,0,0}$z$}%
}}}}
\put(613,-416){\makebox(0,0)[lb]{\smash{{\SetFigFont{11}{13.2}{\rmdefault}{\mddefault}{\updefault}{\color[rgb]{0,0,0}$\delta m$}%
}}}}
\put(778,102){\makebox(0,0)[lb]{\smash{{\SetFigFont{11}{13.2}{\rmdefault}{\mddefault}{\updefault}{\color[rgb]{0,0,0}$x$}%
}}}}
\put(850,-202){\makebox(0,0)[lb]{\smash{{\SetFigFont{11}{13.2}{\rmdefault}{\mddefault}{\updefault}{\color[rgb]{0,0,0}$\delta z$}%
}}}}
\put(2682,-117){\makebox(0,0)[lb]{\smash{{\SetFigFont{11}{13.2}{\rmdefault}{\mddefault}{\updefault}{\color[rgb]{0,0,0}$\ga\delta z$}%
}}}}
\end{picture}%

\end{center}

Assume finally that $D+D'<l(\delta)\leq D+2D'$. See the above
picture. In particular $D'$ is finite, $z\in T$ and $D'= d(y,z)$.
Note that $\delta z$ does not belong to $A_\ga$ since $l(\delta)>
D+D'$, and that 
$$
d(\delta z,x)=d(\delta z,z)-d(z,y)-d(y,x)=l(\delta)- D -D' \leq D'\;.
$$ 
Hence $\ga\delta z\in [y,z]$ and $d(\ga\delta z,y)=d(\delta z,x)=
l(\delta)- D -D'$, so that
$$
d(\ga\delta z,z)=d(z,y)- d(\ga\delta z,y)= 
D'-(l(\delta)- D -D')=D+2D'-l(\delta)\;.
$$ 
Let $m$ be the midpoint of the segment $[\ga\delta z,z]$, so that
$d(m,z)=\frac{1}{2}(D+2D'-l(\delta))$. Hence
$$
d(y,m)= d(y,z)-d(z,m)=\frac{1}{2}(l(\delta)-D)\;.
$$ 
But since $m$ belongs to $A_\delta$ and comes after $z$ on $A_\delta$
oriented by the translation direction of $\delta$, we have
\begin{align*}
d(\delta m,x)=d(\delta m,\delta z)+d(\delta z,x) &=
\frac{1}{2}\big(D+2D'-l(\delta)\big)+\big(l(\delta)- D -D'\big)
\\ & =\frac{1}{2}(l(\delta)-D)=d(y,m)\leq D'\;.
\end{align*}
Hence $\ga\delta m$, which is the point of $[y,z]$ at distance
$d(\delta m,x)$ from $y$, is equal to $m$ and $l(\ga\delta)=0$, as
wanted.  
\end{proof}




\bibliographystyle{plain}
\bibliography{references}

\begin{thebibliography}{10}

\bibitem{AB}
R.~Alperin and H.~Bass.
\newblock Length functions of group actions on {$\Lambda$}-trees.
\newblock In {\em Combinatorial Group Theory and Topology}, volume 111 of {\em
  Ann. of Math. Stud.}, pages 265--378. Princeton Univ. Press, 1987.

\bibitem{A}
S.~Alvarez, D.~Filimonov, V.~Kleptsyn, D.~Malicet, C.~Meni\~{n}o Cot\'{o}n,
  A.~Navas, and M.~Triestino.
\newblock Groups with infinitely many ends acting analytically on the circle.
\newblock {\em J. Topol.}, 12(4):1315--1367, 2019.

\bibitem{B}
A.~F. Beardon.
\newblock Pell's equation and two generator free {M}\"{o}bius groups.
\newblock {\em Bull. London Math. Soc.}, 25(6):527--532, 1993.

\bibitem{Bow}
B.~H. Bowditch.
\newblock Markoff triples and quasi-{F}uchsian groups.
\newblock {\em Proc. London Math. Soc.}, 77(3):697--736, 1998.

\bibitem{LF}
J.~W.~S. Cassels.
\newblock {\em Local Fields}.
\newblock Cambridge University Press, 1986.

\bibitem{CH}
Y.~Cornulier and P.~de~la Harpe.
\newblock {\em Metric Geometry of Locally Compact Groups}.
\newblock European Mathematical Society, 2016.

\bibitem{CM}
M.~Culler and J.~W. Morgan.
\newblock Group actions on {$\mathbb{R}$}-trees.
\newblock {\em Proc. London Math. Soc.}, 55(3):571--604, 1987.

\bibitem{CV}
M.~Culler and K.~Vogtmann.
\newblock The boundary of outer space in rank two.
\newblock In {\em Arboreal Group Theory}, volume~19 of {\em Math. Sci. Res.
  Inst. Publ.}, pages 189--230. Springer, 1991.

\bibitem{SL2}
B.~Eick, M.~Kirschmer, and C.~Leedham-Green.
\newblock The constructive membership problem for discrete free subgroups of
  rank 2 of {${\rm SL_2}(\mathbb{R})$}.
\newblock {\em LMS J. Comput. Math.}, 17(1):345--359, 2014.

\bibitem{GL}
D.~Gaboriau and G.~Levitt.
\newblock The rank of actions on {$\mathbb{R}$}-trees.
\newblock {\em Ann. Sci. \'{E}cole Norm. Sup. (4)}, 28(5):549--570, 1995.

\bibitem{G}
J.~Gilman.
\newblock Two-generator discrete subgroups of {${\rm PSL_2}(\mathbb{R})$}.
\newblock {\em Mem. Amer. Math. Soc.}, 117(561), 1995.

\bibitem{GG}
V.~Guirardel and G.~Levitt.
\newblock Deformation spaces of trees.
\newblock {\em Groups Geom. Dyn.}, 1(2):135--181, 2007.

\bibitem{LU}
R.~C. Lyndon and J.~L. Ullman.
\newblock Groups generated by two parabolic linear fractional transformations.
\newblock {\em Canadian J. Math.}, 21:1388--1403, 1969.

\bibitem{MS}
J.~W. Morgan and P.~B. Shalen.
\newblock Valuations, trees, and degenerations of hyperbolic structures. {I}.
\newblock {\em Ann. of Math. (2)}, 120(3):401--476, 1984.

\bibitem{N}
M.~Newman.
\newblock Pairs of matrices generating discrete free groups and free products.
\newblock {\em Michigan Math. J.}, 15:155--160, 1968.

\bibitem{Paulin}
F.~Paulin.
\newblock The {G}romov topology on {$\mathbb{R}$}-trees.
\newblock {\em Topology Appl.}, 32(3):197--221, 1989.

\bibitem{P}
N.~Purzitsky.
\newblock Two-generator discrete free products.
\newblock {\em Math. Z.}, 126:209--223, 1972.

\bibitem{R2}
G.~Rosenberger.
\newblock Fuchssche {G}ruppen, die freies {P}rodukt zweier zyklischer {G}ruppen
  sind, und die {G}leichung {$x^{2}+y^{2}+z^{2}=xyz$}.
\newblock {\em Math. Ann.}, 199:213--227, 1972.

\bibitem{Serre}
J-P. Serre.
\newblock {\em Trees}.
\newblock Springer, 1980.
\newblock Translated by John Stillwell.

\bibitem{UZ}
M.~Urba\'{n}ski and L.~Zamboni.
\newblock On free actions on {$\Lambda$}-trees.
\newblock {\em Math. Proc. Cambridge Philos. Soc.}, 113(3):535--542, 1993.

\end{thebibliography}

\end{document}